\documentclass[11pt]{article}
\begin{document}
\title{{\bf Some explicit constructions of sets with more sums than 
differences}}
\author{{\sc Peter V. Hegarty} \\ {\em Chalmers University of Technology and 
G\"{o}teborg University} \\ {\em 41296 G\"{o}teborg, Sweden} \\ 
hegarty@math.chalmers.se}
\date{January 5, 2007}
\maketitle

\begin{abstract}
We present a variety of new results on finite 
sets $A$ of integers for which the 
sumset $A+A$ is larger than the difference set $A-A$, so-called MSTD 
(more sums than differences) sets. First we 
show that there is, up to affine transformation, 
a unique MSTD subset of {\bf Z} of size 8. Secondly, starting from some 
examples of size 9, we present several new constructions of infinite families
of MSTD sets. Thirdly we show that for every fixed ordered pair of 
non-negative integers 
$(j,k)$, as $n 
\rightarrow \infty$ a positive proportion of the subsets of 
$\{0,1,2,...,n\}$ satisfy $|A+A| = (2n+1) - j$, $|A-A| = (2n+1) - 2k$.
\end{abstract}

\begin{center} {\bf 1. Introduction} \end{center}
If $A \subseteq {\hbox{{\bf Z}}}$ the sumset and difference set of
$A$ are defined, respectively, as
\begin{eqnarray*}
A + A := \{ x \in {\hbox{{\bf Z}}} : x = a_{1} + a_{2} \; {\hbox{for 
some $a_{1},a_{2} \in A$}} \}, \\
A - A := \{ x \in {\hbox{{\bf Z}}} : x = a_{1} - a_{2} \; {\hbox{for some 
$a_{1},a_{2} \in A$}} \}.
\end{eqnarray*}
If $A$ is a finite set with $n$ elements then it is easy to see that, 
a priori, 
\begin{eqnarray*}
2n-1 \leq |A+A| \leq {n(n+1) \over 2}, \\
2n-1 \leq |A-A| \leq n(n-1) + 1.
\end{eqnarray*}
The upper bounds follow simply from the fact that addition is commutative 
whereas subtraction is not, 
and this certainly suggests that $\lq$most' (speaking informally) finite
sets $A$ should have the property that $|A+A| \leq |A-A|$. A precise
result in support of this intuition was proven by Roesler [6] : 
he showed that for any $n > 0$ and $1 \leq k \leq n$, the average value
of the quotient $|A-A|/|A+A|$, as $A$ runs over the $k$-element subsets
of $\{0,1,...,n\}$, lies in the interval $[1,2)$. A 
question of murky origins, but dating back to the 1960s, is whether 
there exist any finite sets $A$ of integers at all such that $|A+A| > |A-A|$.
The question is interesting because of the apparent difficulty
in finding examples of such sets, whereas it is very easy to 
construct sets $A$ with $|A+A| = |A-A|$. Any arithmetic 
progression has this property or, more generally, any set $A$ with the
property that $A = \{x\} - A$ for some $x \in {\hbox{{\bf Z}}}$. Such a 
set is said to be {\em symmetric (about $x/2$)}. 
\\
\\
Following the practice introduced by Nathanson [3], we refer to 
sets having more sums than differences as MSTD sets. Some examples of such sets
appear in the literature from the late 1960s and 
early 1970s. Conway is said to have 
found the example $A_{1} = \{0,2,3,4,7,11,12,14\}$. The example
$A_{2} = \{0,1,2,4,7,8,12,14,15\}$ appears in [2]. Another interesting
example, $A_{3} = \{0,1,2,4,5,9,12,13,14,16,17,21,24,25,26,28,29\}$, appears
in [1]. Note that for $i = 1,2$, $|A_{i} + A_{i}| - |A_{i} - A_{i}| = 1$, 
whereas $|A_{3} + A_{3}| - |A_{3} - A_{3}| = 4$. 
\\
\\
The following two observations are pertinent :
\\
\\
1. The property of being an MSTD set is invariant under linear mappings 
$x \mapsto ux + v$, $u,v \in {\hbox{{\bf Z}}}$, $u \neq 0$. 
The sizes of sum- and 
difference sets
are unchanged by such mappings. In particular, it suffices to consider finite 
subsets of {\bf Z} which have smallest element 0, and such that 
the greatest common divisor of the elements in the set is 1. Such sets
will be called {\em normalised}.  
\\
\\
2. Starting from any MSTD set $A$ we can construct a sequence $A = A_{1}, 
A_{2},...$ of MSTD sets such that the quotients $|A_{t} + A_{t}|/|A_{t} - 
A_{t}|$ become arbitrarily large. Choose an integer $m$ and set
\begin{eqnarray*}
A_{t} := \left\{ \sum_{i=1}^{t} a_{i} m^{i-1} : a_{i} \in A \right\}.
\end{eqnarray*}
If $m$ is sufficiently large, then $|A_{t}\pm A_{t}| = |A \pm A|^{t}$. This 
method of constructing an infinite family of MSTD sets from a single one
will be called the {\em base expansion method}. 
\\
\\
In particular, these observations imply that any MSTD set gives rise to an
infinite family of such sets. It appears to have been an open problem for some 
time to find some other way of constructing an infinite family of MSTD sets.
In [10], Rusza uses probabilistic
arguments to prove the existence of a multitude of MSTD sets. However, this 
still does not provide explicit constructions. Extending an observation of the 
author and Roesler, such constructions were eventually provided by 
Nathanson [3]. The idea for his type of construction comes from examples like
$A_{1}$ above. Note that that set is the union of a symmetric set 
$\{0,2,3,7,11,12,14\}$ and the single number $4$. The symmetric part contains
an arithmetic progression $\{3,7,11\}$, with some extra numbers tagged on at 
both ends. Nathanson's sets have this type of structure. One starts with a
(proper generalised) arithmetic progression, adds on some structure at both
ends while retaining symmetry, then adds one
further element which results in the sumset being enlarged by one element 
while the 
difference set is left unchanged. 
\par Nathanson's paper also uses a probabilistic method (inspired by Tao) 
to prove the
existence of many MSTD sets in finite abelian groups of the form 
{\bf Z}$/n{\hbox{{\bf Z}}} \times {\hbox{{\bf Z}}}/2{\hbox{{\bf Z}}}$, and 
presents a general method for transforming MSTD sets in finite 
abelian groups to MSTD sets in {\bf Z}. Even more recently, Martin and 
O'Bryant [5]
also use probabilistic methods to prove the following impressive 
result : 
there is a positive 
constant $c$ such that, for all $n >> 0$, at least $c \cdot 2^{n+1}$ subsets 
of $\{0,1,2,...,n\}$ are MSTD sets. 
\\
\\
In summary, probabilistic methods have shown that the phenomenon of MSTD
sets is actually quite common. But explicit constructions of such sets
remain hard to come by, with Nathanson's idea being essentially the only one
in print. This issue of explicitly constructing MSTD sets will be our 
primary concern here. We were motivated by one of the questions Nathanson posed
in his talk [4], namely whether there existed any MSTD sets in {\bf Z} 
of smaller
cardinality than $A_{1}$ above. Our first result is 
\\
\\
{\bf Theorem 1} {\em There are no MSTD subsets of {\bf Z} of size seven.
Moreover, up to linear transformations, $A_{1}$ is the unique such set of size 
8.}
\\
\\
Clearly, the classification of all MSTD sets of a given size, up to 
linear transformation, is a finite computation. To reduce the complexity of the
computation to something manageable, even for very small sizes, is quite
another matter. Our method accomplishes this for sets of size 8, but even
for size 9, the computation was not feasible. We did manage to find all
MSTD sets $A$ of size 9 with the following property : for some element $x$ 
of the sumset $A+A$ there are at least four ordered pairs $(a_{1},a_{2})
\in A \times A$ with $a_{1} + a_{2} = x$. There are exactly nine such 
sets up to linear transformation, namely
\begin{eqnarray*}
A_{4} = \{0,1,2,4,5,9,12,13,14\}, \\
A_{2} = \{0,1,2,4,7,8,12,14,15\}, \\
A_{5} = \{0,2,3,4,7,9,13,14,16\}, \\
A_{6} = \{0,2,3,4,7,11,12,14,16\}, \\
A_{7} = \{0,2,3,4,7,11,15,16,18\}, \\
A_{8} = \{0,2,4,8,9,10,15,17,19\}, \\
A_{9} = \{0,4,6,7,8,14,15,17,21\}, \\
A_{10} = \{0,5,6,9,10,13,16,17,22\}, \\
A_{11} = \{0,4,6,8,11,14,19,21,25\}.
\end{eqnarray*}
Note that $A_{2}$ is the same set as written earlier with the same name. 
$A_{4}$ is a subset of $A_{3}$. The sets $A_{4}, A_{5}, A_{7}$ appear in [3].
Of these, $A_{4}$ and $A_{7}$ are,
along with $A_{1}$, among the infinite family of MSTD sets constructed in 
Theorem 1 of that paper, whereas $A_{5}$ is among the family of MSTD sets
described in Theorem 2 there. The remaining sets $A_{6},A_{8}, 
A_{9},A_{10},A_{11}$ appear to be new. Working from these examples we
will present four constructions of infinite families of MSTD sets
(Theorems 2/3,4,5,6 below) which generalise respectively $A_{11}/A_{9}$, 
$A_{8}$, $A_{6}$ and $A_{4}$. All constructions share common features with,
but are nevertheless different from in a non-trivial sense, those in 
[3] and one another. This reflects the main theme of the paper, namely that
while the most easily describable MSTD sets all seem to have a common core
of features, within this framework there is substantial room for variety. 
The last part of the paper deals with the question 
\\
\par $\lq$How much larger can the sumset be than the difference set ?'
\\
\\
On the one hand, we will answer a question in [5] (a weaker 
version was posed in [3]) by showing (Theorem 8) 
that for every pair $j,k$ of non-negative
integers there is a positive constant $c_{j,k}$ such that, for all 
$n >>_{j,k} 0$, at least $c_{j,k} \cdot 2^{n+1}$ of the subsets $A$ of 
$\{0,1,...,n\}$ satisfy $|A+A| = (2n+1)-j$, $|A-A| = (2n+1)-2k$. 
The proof, which also provides explicit examples for every $j$ and $k$, 
involves two ideas :
firstly, extending the  
probabilistic method in [5], and secondly relating MSTD sets in {\bf Z} to 
MSTD sets in suitably chosen finite cyclic groups. 
\par On the other hand, the base expansion method suggests that the right 
quantity to look at when studying the above question is not
$|A+A| - |A-A|$ but rather 
\begin{eqnarray}
f(A) := {\ln |A+A| \over \ln |A-A|}.
\end{eqnarray}
Results due to 
Freiman-Pigarev [1] and Rusza [7,8,9], establish that
\begin{eqnarray*}
{3 \over 4} \leq f(A) \leq {4 \over 3},
\end{eqnarray*}
for any finite set $A \subseteq {\hbox{{\bf Z}}}$. It is not known if
either bound is sharp, and the state of knowledge is far worse for 
the upper bound. Rusza's
probabilistic method [10] produces a constant $c > 1$ such that there are
a $\lq$multitude' of sets $A$ with $f(A) > c$. He doesn't compute $c$ 
explicitly, but a quick analysis of his method shows that it gives 
$c \approx 1 + 10^{-9}$. The set $A_{3}$
above satisfies $f(A_{3}) = {\ln 59 \over \ln 55} = 1.0175...$ and 
there appears to be nothing in print which beats this. We 
will give explicit examples of sets which do so, if only slightly. 
\\
\\
The remainder of the paper is organised as follows. In Section 2 we will 
prove Theorem 1. 
A mathematica code was written to perform the 
computations necessary to complete the proof.
In Section 3 we will first 
indicate the computations performed 
which allowed us
to conclude that the list of sets $A_{2},...,A_{11}$ was complete in the 
sense mentioned above. Then we will prove Theorems 2-6. In Section 4 we will 
prove Theorem 8 
and exhibit the sets $A$ with larger values of $f(A)$ than 
anything previously written down. In Section 5, we will give some concluding 
remarks and suggestions for further investigations. 

\begin{center} {\bf 2. Proof of Theorem 1} \end{center}
Clearly, it suffices to prove the second statement of the theorem. The 
location of all MSTD subsets of {\bf Z} of a certain size can be 
represented as a finite computation as follows :
\\
\\
Let $A$ be a set of size $n > 0$, with $A = \{a_{i} : i = 1,...,n\}$ and 
$a_{1} > a_{2} > \cdots > a_{n} = 0$. For $i = 1,...,n-1$, represent the
difference $a_{i} - a_{i+1}$ as $\vec{e}_{i}$, the $i$:th standard basis
vector{\footnote{Actually, in our computer program we represent 
$a_{i}-a_{i+1}$ initially by $\vec{e}_{n-i}$.}} 
in {\bf R}$^{n-1}$. A computer program, if now asked to compute those
quantities, will return $|A+A| = n(n+1)/2$, $|A-A| = n(n-1)+1$. So if 
$A$ is to be an MSTD set, there must be a non-trivial coincidence of 
differences. That is, there must exist $i,j,k,l$ such that
$a_{i} - a_{j} = a_{k} - a_{l}$ and $(i,j) \neq (k,l)$. Given such an 
equation we can, by 
projection onto the orthogonal complement in {\bf R}$^{n-1}$ of the 
subspace spanned by $(\vec{e}_{i} - \vec{e}_{j}) - (\vec{e}_{k} - 
\vec{e}_{l})$, now represent the elements of $A$ by vectors in 
{\bf R}$^{n-2}$ and recompute $|A+A|$ and $|A-A|$. If still $|A+A| 
\leq |A-A|$ we can pick another non-trivial identification of elements
in $A-A$, and repeat the above procedure, with the elements of $A$ now 
represented by vectors in {\bf R}$^{n-3}$. Clearly, the computation will
terminate with all MSTD sets of size $n$ and smallest element zero located, 
possibly including infinite parameterised families of such sets. 
\\
\\
The above computation seems to be practically rather unfeasible even 
for $n = 8$ however. We estimate that our machine
would have taken several weeks at least to finish the calculation. 
For $n = 8$ one 
starts with vectors in {\bf R}$^{7}$. It turns out however that, with a 
modest amount of argument, one can show that if $A$ is an MSTD set, then 
there must appear one of 18 possible configurations, each of which
reduces the problem to {\bf R}$^{4}$ or {\bf R}$^{5}$. We thus allowed
our program to instead examine each of these 18 possibilities in turn, and the 
average running time was about 45 minutes. Though the argument used to 
simplify matters could be pushed further, it wasn't obvious to us how to 
do so without an effort which would essentially
balance out the resulting reduction in computing time. Thus our proof of
Theorem 1 will consist of two parts :
\\
\\
{\sc Part One} : reduction to 18 possible cases as described above.
\\
{\sc Part Two} : a computer program to search through all these cases in turn. 

\begin{center} {\bf Part One} \end{center}
As before, let $A \subset {\hbox{{\bf Z}}}$ be a set of size 8, 
$A = \{ a_{i} : i = 1,...,8\}$ where $a_{1} > a_{2} > \cdots > a_{8} = 0$. We
introduce some further notation. Let 
\begin{eqnarray*}
S := \{ (a_{i},a_{j}) : 1 \leq i \leq j \leq 8\}, \\
{\cal S} := A+A = \{a_{i} + a_{j} : (a_{i},a_{j}) \in S\}, \\
D := \{(a_{i},a_{j}) : 1 \leq i < j \leq 8\}, \\
{\cal D} := (A-A)_{> 0} = \{a_{i} - a_{j} : (a_{i},a_{j}) \in D\}.
\end{eqnarray*}
Let ${\cal R}_{1}$ and ${\cal R}_{2}$ be the equivalence relations on $S$ and
$D$ respectively defined by 
\begin{eqnarray*}
\{(a_{i},a_{j}),(a_{k},a_{l})\} \in {\cal R}_{1} \Leftrightarrow a_{i} + a_{j}
= a_{k} + a_{l}, \\
\{(a_{i},a_{j}),(a_{k},a_{l})\} \in {\cal R}_{2} \Leftrightarrow a_{i} - a_{j}
= a_{k} - a_{l}.
\end{eqnarray*}
Let $\Sigma := S/{\cal R}_{1}$ and $\Delta := D/{\cal R}_{2}$. Obviously one
can identify $\Sigma$ with $\cal S$ and $\Delta$ with $\cal D$, but it will
be convenient for us to have a separate notation when referring to the 
equivalence relations. Let $|\Sigma| := \sigma$ and $|\Delta| := \delta$.
Thus $|A+A| = \sigma$ and $|A-A| = 2\delta + 1$. So if $A$ is an MSTD set 
then 
\begin{equation}
2\delta + 1 < \sigma.
\end{equation} 
{\sc Definition :} An equivalence class in $\Sigma$ will be said to be {\em
nice} if it contains at least three elements  
$(a_{i},a_{j})$, $(a_{k},a_{l})$, $(a_{m},a_{n})$ such that all six indices
$i,j,k,l,m,n$ are distinct. 
\\
\\
One readily checks that if there is a nice $\Sigma$-class, then $A$ must 
contain, up to symmetry, one of the following 16 configurations :
\begin{eqnarray}
a_{1} + a_{8} = a_{2} + a_{7} = a_{3} + a_{6}, \\
a_{1} + a_{8} = a_{2} + a_{7} = a_{3} + a_{5}, \\
a_{1} + a_{8} = a_{2} + a_{7} = a_{3} + a_{4}, \\
a_{1} + a_{8} = a_{2} + a_{7} = a_{4} + a_{5}, \\
a_{1} + a_{8} = a_{2} + a_{6} = a_{3} + a_{5}, \\
a_{1} + a_{8} = a_{2} + a_{6} = a_{3} + a_{4}, \\
a_{1} + a_{8} = a_{2} + a_{6} = a_{4} + a_{5}, \\
a_{1} + a_{8} = a_{2} + a_{5} = a_{3} + a_{4}, \\
a_{1} + a_{8} = a_{3} + a_{6} = a_{4} + a_{5}, \\
a_{1} + a_{7} = a_{2} + a_{6} = a_{3} + a_{5}, \\
a_{1} + a_{7} = a_{2} + a_{6} = a_{3} + a_{4}, \\
a_{1} + a_{7} = a_{2} + a_{6} = a_{4} + a_{5}, \\
a_{1} + a_{7} = a_{2} + a_{5} = a_{3} + a_{4}, \\
a_{1} + a_{7} = a_{3} + a_{6} = a_{4} + a_{5}, \\
a_{1} + a_{6} = a_{2} + a_{5} = a_{3} + a_{4}, \\
a_{2} + a_{7} = a_{3} + a_{6} = a_{4} + a_{5}.
\end{eqnarray}
Let $n_{1} \leq n_{2} \leq \cdots \leq n_{\delta}$ be the sizes of the 
$\Delta$-classes, arranged in some increasing order. Hence 
\begin{equation}
\sum_{i = 1}^{\delta} n_{i} = |D| = 28.
\end{equation}
Note that $n_{1} = 1$ since $(a_{1},a_{8})$ is in a class by itself. If 
$(a_{i},a_{j})$ is ${\cal R}_{2}$-equivalent to 
$(a_{k},a_{l})$, where $i \leq k \leq j \leq l$, then 
$(a_{i},a_{l})$ is ${\cal R}_{1}$-equivalent to $(a_{k},a_{j})$. This defines a 
mapping $\phi$ from ${\cal R}_{2}$-equivalent pairs of elements of $D$ 
to ${\cal R}_{1}$-equivalent pairs of elements of $S$. The mapping is obviously
at most 2-1, hence its range consists of at least $\frac{1}{2} 
\sum_{i = 1}^{\delta} \left( \begin{array}{c} n_{i} \\ 2 \end{array} \right)$
pairs. If, in addition, there are no nice $\Sigma$-classes, it follows 
that 
\begin{equation}
\sigma \leq |S| - {1 \over 2} \sum_{i = 1}^{\delta} \left( \begin{array}{c}
n_{i} \\ 2 \end{array} \right) = 36 - {1 \over 2} \sum_{i=1}^{\delta}
\left( \begin{array}{c} n_{i} \\ 2 \end{array} \right).
\end{equation}
We examine the following two cases :
\\
\\
{\sc Case I} : $n_{2} > 1$.
\\
\\
{\sc Case II} : $n_{2} = 1$ and there are no nice $\Sigma$-classes. 
\\
\\
First consider {\sc Case I}. 
Then every $\Delta$-class, other than that consisting of 
the single pair $(a_{1},a_{8})$, contains at least two members of $D$. First of
all this forces $a_{1} - a_{7} = a_{2} - a_{8}$, hence 
\begin{equation}
a_{1} + a_{8} = a_{2} + a_{7}.
\end{equation}
Next consider the three differences $a_{1} - a_{6}$, $a_{2} - a_{7}$ and 
$a_{3} - a_{8}$. The two largest ones must be equal. If 
$a_{1} - a_{6} = a_{3} - a_{8}$ then $a_{1} + a_{8} = a_{3} + a_{6}$ which, 
together with (21), implies that $A$ contains the nice configuration (3). 
Otherwise, up to symmetry, we may assume that $a_{1} - a_{6} = a_{2} - a_{7}
> a_{3} - a_{8}$. It is then easily checked that one of the following 
three possibilities must hold :
\\
\\
(a) $a_{1} - a_{5} = a_{3} - a_{8}$. This, together with (21), implies that 
$A$ contains the nice configuration (4). 
\\
\\
(b) $a_{1} - a_{5} = a_{2} - a_{6}$. Thus in $A$ we have that 
\begin{equation}
a_{1}-a_{2} = a_{5}-a_{6} = a_{6}-a_{7}=a_{7}-a_{8}.
\end{equation}
(c) $a_{2}-a_{6} = a_{3} - a_{8}$. Thus in $A$ we have that 
\begin{equation}
a_{1}-a_{2}=a_{6}-a_{7}=a_{7}-a_{8}=\frac{1}{2} (a_{2}-a_{3}).
\end{equation}
This completes the analysis of {\sc Case I}. We have shown that, under these
circumstances, either there is a nice $\Sigma$-class or, up to 
symmetry, $A$ contains one of
the configurations (22) and (23). 
\\
\\
Next consider {\sc Case II}. We claim that in this case, $A$ cannot be an MSTD 
set. Suppose otherwise. A priori, $\sigma \leq 36$ so (2) forces $\delta 
\leq 17$. Now this plus (19) and the assumption that $n_{2} = 1$ mean that
$\sum_{i = 1}^{\delta} \left( \begin{array}{c} n_{i} \\ 2 \end{array} \right)
\geq 11$. Then (20) implies that $\sigma \leq 30$. The idea is now to iterate
this kind of argument to gradually reduce $\delta$ until we obtain the 
contradiction that $\delta < 7$.
\par If $\sigma \leq 30$ then (2) forces $\delta \leq 14$. Then (19) and 
the assumption that $n_{2} = 1$ force $\sum_{i = 1}^{\delta}
\left( \begin{array}{c} n_{i} \\ 2 \end{array} \right) \geq 16$. Thus (20)
forces $\sigma \leq 28$. 
\par If $\sigma \leq 28$ then (2) forces $\delta \leq 13$. Then (19) and 
the assumption that $n_{2} = 1$ force $\sum_{i = 1}^{\delta}
\left( \begin{array}{c} n_{i} \\ 2 \end{array} \right) \geq 19$. Thus (20)
forces $\sigma \leq 26$.
\par If $\sigma \leq 26$ then (2) forces $\delta \leq 12$. Then (19) and 
the assumption that $n_{2} = 1$ force $\sum_{i = 1}^{\delta}
\left( \begin{array}{c} n_{i} \\ 2 \end{array} \right) \geq 22$. Thus (20)
forces $\sigma \leq 25$.
\par If $\sigma \leq 25$ then (2) forces $\delta \leq 11$. Then (19) and 
the assumption that $n_{2} = 1$ force $\sum_{i = 1}^{\delta}
\left( \begin{array}{c} n_{i} \\ 2 \end{array} \right) \geq 25$. Thus (20)
forces $\sigma \leq 23$.
\par If $\sigma \leq 23$ then (2) forces $\delta \leq 10$. Then (19) and 
the assumption that $n_{2} = 1$ force $\sum_{i = 1}^{\delta}
\left( \begin{array}{c} n_{i} \\ 2 \end{array} \right) \geq 30$. Thus (20)
forces $\sigma \leq 21$.
\par If $\sigma \leq 21$ then (2) forces $\delta \leq 9$. Then (19) and 
the assumption that $n_{2} = 1$ force $\sum_{i = 1}^{\delta}
\left( \begin{array}{c} n_{i} \\ 2 \end{array} \right) \geq 36$. Thus (20)
forces $\sigma \leq 18$.
\par If $\sigma \leq 18$ then (2) forces $\delta \leq 8$. Then (19) and 
the assumption that $n_{2} = 1$ force $\sum_{i = 1}^{\delta}
\left( \begin{array}{c} n_{i} \\ 2 \end{array} \right) \geq 44$. Thus (20)
forces $\sigma \leq 14$.
\par If $\sigma \leq 14$ then (2) forces $\delta \leq 6$ and we have our 
desired contradiction. 
\\
\\
This completes the analysis of {\sc Case II}. We have thus shown that, if $A$
is an MSTD set then, up to symmetry, 
either it contains one of the 16 nice configurations
(3)-(18), or one of the configurations (22) and (23). 
Thus the proof of Theorem 1
is reduced to 18 possible cases, as claimed, and we have completed Part 
One of the proof. 

\begin{center} {\bf Part Two} \end{center}
If $A$ contains a nice configuration, then the differences 
$a_{i} - a_{i+1}$ can be represented with vectors in {\bf R}$^{5}$. If either 
(22) or (23) hold, then it suffices with {\bf R}$^{4}$. We can then
write a program to search for normalised MSTD sets as outlined above. The 
rather unwieldy 
mathematica code for such a program, with each of the 18 possible
input configurations, can be obtained from the author. 
The total running time for
the program on our network was about 15 hours. The only outputted MSTD sets
were $A_{1}$ and $14 - A_{1}$, thus proving Theorem 1.

\begin{center} {\bf 3. Explicit families of MSTD sets} \end{center}
If $A$ is an MSTD set of size 9, then the same type of reasoning as above can 
be used to show that either there is a nice $\Sigma$-class or $A$
contains one of a small number of other possible configurations analogous to
(22) and (23). However, a nice $\Sigma$-class only allows us now to represent
the differences $a_{i} - a_{i+1}$, a priori, as vectors in {\bf R}$^{6}$. 
Plus, the number of possible configurations is now much larger. The resulting
computation was not practically feasible with our code, so we have not obtained
with certainty a full classification of all MSTD sets of size 9, up to linear 
transformations. Instead, we classified all such sets in which some 
$\Sigma$-class contains at least 4 elements. One readily verifies that
this leaves, up to symmetry, 25 possible configurations in $A$ analogous to
those in eqs.(3)-(18). 
The computer then located, over a period of about two weeks and after
multiple crashes, the nine normalised MSTD sets $A_{2},A_{4},A_{5},...$
listed in the introduction. We now present a sequence of constructions 
of infinite families of MSTD sets which include five of these nine examples :
\\
\\
{\bf Theorem 2} {\em Let $n \geq 3$, $0 \leq k \leq n-2$. Let $d > 1$ be a
divisor of $2^{n-k} - 1$. Set 
\begin{eqnarray*}
X := \{ 2^{n} - 2^{j} : k \leq j \leq n\}, \\
m := (2^{n+1}+d) + (2^{n} - 2^{k}), \\
Y := m - X, \\
Z := \left\{2^{n} + jd : 1 \leq j \leq \frac{2^{n}-2^{k}}{d} \right\}, \\
B := X \sqcup Y \sqcup Z, \\
a := 2^{n}, \\
A := B \cup \{a\}.
\end{eqnarray*}
Then $A$ is a normalised MSTD set with $|A+A| = |A-A| + 1$.}
\\
\\
{\sc Proof} : One may verify the following facts :
\par (i) $B$ is a symmetric set, $B = m-B$.
\par (ii) $A-A = B-B$.
\par (iii) $A+A = (B+B) \sqcup \{2a\}$.
\\
The verifications involve calculations similar to those appearing in 
the proofs of Theorems 1-4 in [3], hence are omitted. 
\\
\\
{\sc Remark} : The set $A_{11}$ is the case $n=3,k=1,d=3$ of the above theorem.
\\
\\
In the notation of Theorem 2, when $k = 0$ there is a simpler, similar
construction of a family of normalised 
MSTD sets : basically the arithmetic progression $Z$ is not needed.
\\
\\
{\bf Theorem 3} {\em Let $n \geq 3$ and $1 \leq l \leq n-2$. Set 
\begin{eqnarray*}                     
X := \{ 2^{n} - 2^{j} : 0 \leq j \leq n\}, \\
m := (2^{n+1}-1) + (2^{n} - 2^{l}), \\
Y := m-X, \\
B := X \sqcup Y, \\
a := 2^{n}, \\
A := B \sqcup \{a\}.
\end{eqnarray*}
Then $A$ is a normalised MSTD set with $|A+A| = |A-A| + 1$}.
\\
\\
{\sc Proof} : As for Theorem 2. Observe that $A_{9}$ is the case $n=3, l=1$ of
the theorem.
\\
\\
Next we turn to a generalisation of $A_{8}$ :
\\
\\
{\bf Theorem 4} {\em Let $n,d > 1$. Set 
\begin{eqnarray*}
X := \{jd : 0 \leq j \leq n\}, \\
m := (4n+1)d + 1, \\
Y := m-X, \\
Z := (2nd,(2n+1)d], \\
B := X \sqcup Y \sqcup Z, \\
a := 2nd, \\
A := B \sqcup \{a\}.
\end{eqnarray*}
Then $A$ is a normalised MSTD set with $|A+A| = |A-A| + 1$.}
\\
\\
{\sc Proof} : One may verify that $B$ is symmetric, that 
$A-A = B-B$ and that $A+A = (B+B) \sqcup \{2a\}$.
\\
\\
{\sc Remark} : The set $A_{8}$ is the case $n=d=2$ of this theorem. 
\\
\\
Next we turn to $A_{6}$. Note that $A_{6} = A_{1} \cup \{16\}$ and that it does
not have a symmetric subset of size 8. The following generalises $A_{6}$ 
in a different direction to the generalisation of $A_{1}$ given in [3] :
\\
\\
{\bf Theorem 5} {\em Let $n \geq 2$, $k \geq 3$. Set 
\begin{eqnarray*}
X := \{2j : 0 \leq j < n\}, \\
m := 2(k+1)n - 2, \\
Y := m - X, \\
Z := \{2jn - 1 : 1 \leq j \leq k\}, \\
B := X \sqcup Y \sqcup Z, \\
a^{*} := 2n, \\
A^{*} := B \sqcup \{a^{*}\}, \\
a := m+2, \\
A := A^{*} \sqcup \{a\}.
\end{eqnarray*}
Then $A^{*}$ and $A$ are both MSTD sets and $|A^{*}+A^{*}| - 
|A^{*}-A^{*}| = |A+A| - |A-A| = 1$.}
\\
\\
{\sc Remark} : The sets $A$ in Theorem 5 
provide explicit examples of 
MSTD sets which are not obtained from a symmetric subset by adding a single
element, followed possibly by a base expansion. 
Indeed $B$ is a maximal symmetric subset of $A$ and $|A\backslash B|
= 2$. Finally note that the set $A_{6}$ is the case $n=2,k=3$ of the 
theorem.
\\
\\
{\sc Proof of Theorem 5} : One may verify the following facts, from which
the theorem follows :
\\
\\
(i) $B$ is symmetric.
\\
(ii) $A^{*}-A^{*} = B-B$.
\\
(iii) $A^{*}+A^{*} = (B+B) \sqcup \{2a^{*}\}$.
\\   
(iv) $A-A = (A^{*}-A^{*}) \sqcup \{\pm a, \pm [a-(a^{*}-1)] \}$. 
\\
(v) $A+A = (A^{*}+A^{*}) \sqcup \{2a, 2a-2, a+(2kn-1), a+a^{*} \}$.
\\
\\
Our fourth construction is a novel generalisation of 
the set $A_{4}$, different from that covered by Theorem 1 of [3] :
\\
\\
{\bf Theorem 6} {\em Let $n,k > 1$. Set 
\begin{eqnarray*}
X := [0,n], \\
m := (2k+3)n, \\
Y := m-X, \\
Z := \sqcup_{j=1}^{k} [2jn+1,(2j+1)n-1], \\
B := X \sqcup Y \sqcup Z, \\
a := 2n, \\
A := B \sqcup \{a\}.
\end{eqnarray*}
Then $A$ is an MSTD set with $|A+A| - |A-A| = 1$. Furthermore, set 
\begin{eqnarray*}
W := [(2k+4)n,(2k+5)n), \\
{\cal A} := A \sqcup W.
\end{eqnarray*}
Then $\cal A$ is also an MSTD set and $|{\cal A} + {\cal A}| - |{\cal A}
- {\cal A}| = 2$.}
\\
\\
{\sc Proof} : One may verify the following facts, from which the theorem 
\\ follows :
\\
\\
(i) $B$ is symmetric.
\\
(ii) $A-A = B-B$.
\\
(iii) $A+A = (B+B) \sqcup \{2a\}$.
\\
(iv) ${\cal A} - {\cal A} = (A-A) \sqcup \pm [(2k+3)n+1,(2k+5)n-1]$. Hence
$|{\cal A} - {\cal A}| = |A-A| + (4n-2)$.
\\
(v) ${\cal A} + {\cal A} = (A+A) \sqcup [(4k+6)n+1,(4k+10)n-2]$. Hence
$|{\cal A}+{\cal A}| = |A+A| + (4n-1)$.
\\
\\
{\sc Remarks} : (a) In the case $n=k=2$, the set $A$ coincides with 
$A_{4}$ and the set $\cal A$ is $A_{4} \cup \{16,17\} = 
\{0,1,2,4,5,9,12,13,14,16,17\}$. Denote this set by $A_{12}$ for future 
reference.
\\
(b) The sets $\cal A$ provide an explicit infinite
family of examples of 
MSTD sets in which the size of the sumset is at least two more than the 
size of the difference set, and 
which are not obtained from a fixed MSTD set by either a 
base expansion or by the method of Lemma 7 below (the latter method is implicit
in [3]).  
\\
(c) The set $A_{2}$ has as a superset $A_{13} = A_{2} \cup \{18,19,20\}$ which
also satisfies $|A_{13}+A_{13}| - |A_{13}-A_{13}| = 2$. This set is not, 
however, covered by Theorem 6, and we have not spotted a way to 
generalise either $A_{13}$ or $A_{2}$ itself.
\\ 
(d) Given an MSTD set $M$, let $s(M)$ denote the maximal size of a symmetric
subset of $M$. In Theorem 5, we gave examples in which $|M| - s(M) = 2$. 
All our examples prior to that satisfy $|M| - s(M) = 1$. 
\par Now consider the sets $\cal A$ in Theorem 6. Each $\cal A$ 
has the symmetric subset $B$ and $|{\cal A} \backslash B| =
n+1$. Let $C = \{n,(2k+3)n\} \sqcup (Z \backslash [2n+1,3n-1])$. Then
${\cal A} \backslash C$ is also symmetric and $|{\cal A} \backslash C|
= (k-1)(n-1) + 2$. It is easily checked that $\cal A$ has no symmetric 
subset which is larger than both $B$ and ${\cal A} \backslash C$. 
Thus we get examples of MSTD sets $M$ in which the difference 
$|M| - s(M)$ can be made arbitrarily large. The sets constructed in 
the proof of Theorem 8 below provide further examples of this 
phenomenon. Again, the point is that in both cases, these 
examples are not base expansions of a fixed set.

\begin{center} {\bf 4. How much larger can the sumset be ?} \end{center}
Starting from any fixed MSTD set, the base expansion method allows one 
to construct MSTD sets $A$ for which the quotients $|A+A|/|A-A|$ become
arbitrarily large. Here we are interested in both the quotient 
$|A+A|/|A-A|$ and the difference $|A+A| - |A-A|$. Our first step is to
describe a way to identify an  
MSTD subset $A$ of {\bf Z}
with MSTD subsets of {\bf Z}$/n{\hbox{{\bf Z}}}$ for suitable $n$ depending 
on $A$. This idea is also implicitly contained in [3]. 
\par Let $A$ be a normalised MSTD subset of {\bf Z}, $m := 
\max \{a \in A\}$. Then $A$ can be 
considered as an MSTD subset of {\bf Z}$/n{\hbox{{\bf Z}}}$ for any $n > 2m$.
More generally, we make the following definition :
\\
\\
{\sc Definition} : Let $A$ be a normalised subset of {\bf Z} (not necessarily
an MSTD set) with 
largest element $m$, and $n > 0$ an integer.
Set $A^{\prime} := A \cap [0,n)$ and identify $A^{\prime}$ with a subset of
{\bf Z}$/n{\hbox{{\bf Z}}}$ in the natural way. We call $A^{\prime}$ the 
{\em reduction} of $A$ modulo $n$.  
Then $A$ is said to be {\em reducible} modulo $n$ if 
$A = \{ x \in [0,m] : x \; ({\hbox{mod $n$}}) \in A^{\prime} \}$. If, in 
addition, $A^{\prime}$ is an MSTD set in {\bf Z}$/n{\hbox{{\bf Z}}}$, we say 
that $A$ has {\em good MSTD-reduction} modulo $n$.
\\
\\
The important observation, which is basically a reformulation and 
sharpening of Theorem 8 of [3] in the case of cyclic groups, is the following :
\\
\\
{\bf Lemma 7} {\em Let $A$ be a normalised subset of {\bf Z} and $n > 0$ an 
integer. Suppose $A$ is reducible modulo $n$ with reduction $A^{\prime}$.
Let 
\begin{eqnarray*}
f(n,A) := |A^{\prime} + A^{\prime}| - |A^{\prime}-A^{\prime}|.
\end{eqnarray*}
Put $B := A \cup (A+n)$. Then $B$ is a normalised subset of {\bf Z} and 
\begin{equation}
|B+B| - |B-B| = (|A+A| - |A-A|) + 2 \cdot f(n,A).
\end{equation}}
{\sc Proof} : Let $X$ be any congruence class of integers modulo $n$. The
reducibility of $A$ modulo $n$ implies the following : If $x_{1},x_{2} \in 
X \cap A+A$ (resp. $X \cap A-A$), and $x_{3} \in (x_{1},x_{2}) \cap X$, then 
$x_{3} \in X \cap A+A$ (resp. $X \cap A-A$). Furthermore, $X \cap A+A$ 
(resp. $X \cap A-A$) is non-empty if and only if $X \cap B+B$ (resp. 
$X \cap B-B$) is. 
\par Suppose $X \cap A+A$ is non-empty with largest element $x$. It follows
that $X \cap B+B = (X \cap A+A) \sqcup \{x+n, x+2n\}$.
Similarly, suppose $X \cap A-A$ is non-empty with largest element
$x_{0}$ and smallest element $x_{1}$. Then $X \cap B-B = (X \cap A-A)
\sqcup \{x_{0}+n,x_{1}-n\}$.  
\par Eq. (24) follows immediately from these observations. 
\\
\\
Each of the sets listed on page 4 has good MSTD reduction modulo $n$, for 
some $n$ considerably less than $2 \cdot \max \{a \in A\}$. For example, we
can take 
\begin{eqnarray*}
A_{4}^{\prime} = A_{4} \backslash \{12,13,14\}, 
\;\; n=12, \;\; f(12,A_{4}) = 1, \\
A_{2}^{\prime} = A_{2}, \;\; n =18, \;\; f(18,A_{2}) = 1, \\
A_{5}^{\prime} = A_{5} \backslash \{16\}, \;\; n = 16, \;\; f(16,A_{5}) = 1, \\
A_{6}^{\prime} = A_{6}, \;\; n=21, \;\; f(21,A_{6}) = 2, \\
A_{7}^{\prime} = A_{7} \backslash \{16,18\}, \;\; n=16, \;\; f(16,A_{7}) = 1, \\
A_{8}^{\prime} = A_{8} \backslash \{17,19\}, \;\; n=17, \;\; f(17,A_{8}) = 1, \\
A_{9}^{\prime} = A_{9} \backslash \{17,21\}, \;\; n=17, \;\; f(17,A_{9}) = 1, \\
A_{10}^{\prime} = A_{10} \backslash \{17,22\}, \;\; n=17, \;\; f(17,A_{10}) = 1, 
\\
A_{11}^{\prime} = A_{11} \backslash \{21,25\}, \;\; n=21, \;\; f(21,A_{11}) = 1.
\end{eqnarray*}
More generally, in the notation of the theorems above :
\\
\\
Theorem 2 : For given $n,k,d$ the set $A$ has good MSTD reduction 
modulo $r := (2^{n+1}+d) + (2^{n-1}-2^{k})$ and $f(r,A) = 1$.
\\
Theorem 3 : For given $n,l$ the set $A$ has good MSTD reduction modulo 
$r := (2^{n+1}-1)+(2^{n-1}-2^{l})$ and $f(r,A) = 1$.
\\
Theorem 4 : For given $n,d$ the set $A$ has good MSTD reduction modulo 
$r := 4nd+1$ and $f(r,A) = 1$. 
\\
Theorem 5 : For given $n,k$ the set $A$ has good MSTD reduction modulo 
$r := 2(k+2)n + 1$. If $k = 3$ then $f(r,A) = 2$. If $k > 3$ then 
$f(r,A) = 1$.
\\
Theorem 6 : For given $n,k$ the set $A$ has good MSTD reduction modulo 
$r := (2k+2)n$ and $f(r,A) = 1$. Note that, in fact, $A+A \; ({\hbox{mod $r$}})
= {\hbox{{\bf Z}}}/r{\hbox{{\bf Z}}}$ and $A-A \; ({\hbox{mod $r$}}) = 
{\hbox{{\bf Z}}}/r{\hbox{{\bf Z}}} \backslash \{r/2\}$.
\\
\\
The rather tedious proofs of these statements are omitted. 
Instead we note that of particular interest is the fact that 
$f(21,A_{6}) > 1$. This means that
in repeatedly 
applying Lemma 7, the difference in size between the sum- and difference
sets will grow more quickly. Let 
\begin{eqnarray*}
A_{14} := A_{6} \sqcup (A_{6} + 21) \sqcup (A_{6} + 42) \end{eqnarray*}
\begin{eqnarray*} 
= 
\{0,2,4,5,9,12,13,14,16,21,23,25,26,30,33,34,35,37,42,44,46,47,51,54,55,56,
58\}.
\end{eqnarray*}
Then $|A_{14}+A_{14}| = 114$ and $|A_{14}-A_{14}| = 105$, so in the notation of 
(1), we have 
$f(A_{14}) = {\ln 114 \over \ln 105} = 1.01767...$. Thus $f(A_{14}) > f(A_{3})$.
By the way, note that $A_{3} = A_{12} \cup (A_{12} + 12)$.
\par The following example does even better : Let 
\begin{eqnarray*}
X := 
\{0,1,2,4,5,9,12,13,17,20,21,22,24,25\}.
\end{eqnarray*}
(this set appears in [5], where it is 
denoted $S_{4}$). Then $|X+X| - |X-X| = 4$, $X$ is reducible modulo $20$ and 
$f(20,X) = 2$. 
Take 
\begin{eqnarray*}
A_{15} := X \cup (X + 20) \end{eqnarray*}
\begin{eqnarray*}
= \{0,1,2,4,5,9,12,13,17,20,21,22,24,25,29,32,33,37,40,41,42,44,45\}.
\end{eqnarray*}
Then $|A_{15}+A_{15}| = 91$ and $|A_{15}-A_{15}| = 83$, so $f(A_{15}) = 
{\ln 91 \over \ln 83} = 1.02082...$.  
\\
\\
Our final result resolves in the affirmative Conjecture 20 of [5] :
\\
\\
{\bf Theorem 8} {\em Let $j,k$ be any two non-negative integers. Then there
exists a constant $c_{j,k} \in (0,1)$ such that, for all $n$ sufficiently
large, depending on $j$ and $k$, at least $c_{j,k} \cdot 2^{n+1}$ of
the subsets $A$ of $\{0,1,2,...,n\}$ satisfy $|A+A| = (2n+1)-j$ and 
$|A-A| = (2n+1)-2k$.}
\\
\\
{\sc Proof} : The proof is an extension of the method of [5] and consists
of two separate constructions depending on the sign of $j-k$. 
\\
\\
{\sc Case I} : $j \leq k$. 
\\
\\
We make use of the set $A_{12} = \{0,1,2,4,5,9,12,13,14,16,17\}$ above. It
satisfies $A_{12} + A_{12} = [0,34]$ and $A_{12} - A_{12} = [-17,17] \backslash 
\{\pm 6\}$. Furthermore, this set has good MSTD reduction modulo 12. 
\par Let $j,k$ be given. Set 
\begin{eqnarray*}
X_{k} := \bigcup_{t=0}^{k-1} A_{12} + 12t.
\end{eqnarray*}
Let $m_{k} = \max X_{k} = 12k+5$. By the proof of Lemma 7, we have that
$X_{k} + X_{k} = [0,2m_{k}]$ and $|X_{k}-X_{k}| = (2m_{k}+1) - 2k$, since the
difference set misses all numbers congruent to 6 (mod 12) in the interval
$[-m_{k},m_{k}]$. Now set
\begin{eqnarray*}
X_{j,k} := X_{k} \backslash \{12(k-t)+1 : 0 \leq t < j \}.
\end{eqnarray*}
One readily checks that $X_{j,k} - X_{j,k} = X_{k}-X_{k}$ but that
$X_{j,k}+X_{j,k} = [0,2m] \backslash \{24(k-t)+3 : 0 \leq t < j\}$. Thus the 
sets $X_{j,k}$ already provide explicit examples of sets satisfying the 
requirements of Theorem 8. To prove the existence of positive 
constants $c_{j,k}$ we use the method of [5]. Set $L_{j,k} := (m_{k} - 
X_{j,k}) \sqcup (m_{k},2m_{k}]$ and $U_{j,k} := 
n - (X_{k} \sqcup (m_{k},2m_{k}])$, for any $n \geq 4m_{k}+1$. Finally put
\begin{eqnarray*}
A_{j,k} := L_{j,k} \sqcup R \sqcup U_{j,k},
\end{eqnarray*}
where $R$ is the random subset of $(2m_{k},n-2m_{k})$ obtained by choosing
each number in the interval independently with probability $1/2$. Our choices
of $L_{j,k}$ and $U_{j,k}$ imply that 
\begin{eqnarray*}
\# \{ ([0,4m_{k}] \cup [2n-4m_{k},2n]) \backslash (A_{j,k}+A_{j,k}) \} = j, \\
\# \{ [n-4m_{k},n] \backslash (A_{j,k}-A_{j,k}) \} = k.
\end{eqnarray*}
It then suffices to apply the same type of argument as in [5] to show that, 
with high probability, both $(4m_{k},2n-4m_{k}) \subset A_{j,k}+A_{j,k}$
and $[0,n-4m_{k}) \subset A_{j,k} - A_{j,k}$, and thus deduce the 
existence of a constant $c_{j,k} > 0$.     
\\
\\
{\sc Case II} : $j \geq k$. 
\\
\\
We start by describing,
for each $j \geq 0$, an integer $m_{j} \geq -1$ and a subset $L_{j} \subset
[0,m_{j}]$ with the following properties :
\\
\\
(i) $L_{j} + L_{j} = [0,2m_{j}] \backslash \{2x_{1},...,2x_{j}\}$ where 
$0 < x_{1} < \cdots < x_{j} < m_{j}/2$.
\\
(ii) None of the numbers $2x_{u} - x_{v}$, for $1 \leq u,v \leq j$, is in 
$L_{j}$. 
\\
\\
Set $m_{0} := -1$, $L_{0} := \phi$. This is consistent with (i) and (ii) 
above. We define the numbers 
$m_{1} < m_{2} < \cdots$ and the sets $L_{1} \subset L_{2} \subset \cdots$
inductively as follows : for each $i \geq 0$, 
\begin{eqnarray*}
m_{i+1} := 3m_{i} + 16,
\end{eqnarray*}
\begin{eqnarray*}
L_{i+1} := L_{i} \sqcup \left\{ (m_{i} + 1) + 
\{0,1,2,5,m_{i}+10,m_{i}+11,...,2m_{i}+15\} \right\}.
\end{eqnarray*}
Each $m_{j}$ is an odd 
number and one readily checks that
\begin{eqnarray}
m_{j} = 7 \cdot 3^{j} - 8, \end{eqnarray}
\begin{eqnarray}
|L_{j}| = {7 \over 2} (3^{j}-1) + 2j,
\end{eqnarray}
\begin{eqnarray}
L_{j} + L_{j} = [0,2m_{j}] \backslash \{2x_{t} = 2(m_{t} + 5) : 0 \leq t < j \}.
\end{eqnarray}
It is then a simple exercise to show that conditions (i) and (ii) above on
$L_{j}$ are satisfied. 
\\
\\
Now let $j,k$ be given with $j \geq k$ and first suppose $k > 0$. Set 
\begin{eqnarray*}
L_{j,k} := L_{j} \sqcup \{2x_{t+1}-x_{t} : 1 \leq t \leq j-k\},
\end{eqnarray*}
and observe that the numbers $x_{t}$ grow sufficiently fast so that
for any $k$, $L_{j,k}+L_{j,k} = L_{j}+L_{j}$. For any $n \geq 2m_{j}+1$ set
\begin{eqnarray*}
U_{j,k} = U_{j} := n - (L_{j} \cup \{x_{1},...,x_{j}\}).
\end{eqnarray*}
Then, as usual, take finally $A_{j,k} := L_{j,k} \sqcup R \sqcup U_{j,k}$,
where $R$ is the
random subset of $(m_{j},n-m_{j})$ obtained by choosing each number 
independently with probability $1/2$. One easily checks that 
\begin{eqnarray*}
([0,2m_{j}] \cup [2n-2m_{j},2n]) \backslash (A_{j,k}+A_{j,k}) = \{2x_{t} : 
t = 1,...,j\}, 
\\ {\hbox{and}} \;\;
[n-2m_{j},n] \backslash (A_{j,k}-A_{j,k}) 
= \{n-1-2x_{t} : t=1 \; {\hbox{or}} \; j-k+2 \leq t \leq j\}.
\end{eqnarray*}
By the method of [5] both
$(2m_{j},2n-2m_{j}) \subset A_{j,k}+A_{j,k}$ and 
$[0,n-2m_{j}) \subset A_{j,k}-A_{j,k}$ occur with high probability, from which
we deduce the existence of a positive constant $c_{j,k}$.
\par When $k = 0$ we just have to be a little careful. One may check that
the following choices work :
\begin{eqnarray*}
L_{j,0} := L_{j+2} \sqcup \{x_{1}\} \sqcup \{2x_{t+1}-x_{t} : 1 \leq t \leq j+1\},
\;\;\;\; U_{j,0} := U_{j+2}, 
\end{eqnarray*} 
and thus the proof of Theorem 8 is complete. 
\\
\\
{\sc Remarks} : (a) In each of the above constructions we get explicit 
examples of sets satisfying the requirements of Theorem 8 by taking the set
$R$ to consist of the entire interval over which it is defined. These 
examples are not covered by any of the constructions in Section 3, and in 
particular that in {\em Case II} is somewhat different from any of the 
earlier ones. 
\\
(b) From the above proof we obtain estimates of the form  
$c_{j,k} = \Omega (2^{- \Theta(m_{l})})$, where $l = \max \{j,k\}$. The 
numbers $m_{l}$ grow linearly when $j < k$ but exponentially when $k > j$.
This lack of symmetry is unsatisfying, though in any case the estimates are
likely to be way smaller 
than the truth for all values of $j$ and $k$.      

\begin{center} {\bf 5. Concluding remarks} \end{center}
We have provided various explicit constructions of non-trivial infinite 
families of MSTD sets. All our constructions, including those in the proof
of Theorem 8, and in common with 
those already in the literature, are roughly based on 
some type of symmetric set which is $\lq$perturbed' slightly by adding on a 
small number of elements. Often, though not always, the symmetric set
is itself constructed out of some (generalised) arithmetic progression. 
Theorems 2,3 provide an example where this is not so : there is an 
arithmetic progression $Z$ in the former construction, but both are
based on a geometric progression $X$. The sets constructed in the proof of 
Theorem 8 are also of a somewhat different character. 
Note that in Theorem 6, the
set $Z$ is a GAP of dimension 2. We have not investigated whether 
the arithmetic progressions appearing in our various constructions 
can be replaced, as in [3], by GAP:s of higher dimension. More 
interesting, though, would be to have explicit examples of MSTD sets 
which are, in some meaningful 
sense, $\lq$radically' different from the blueprint of
a perturbed symmetric set. The ubiquity of MSTD sets, as exhibited by 
probabilistic techniques to which Theorem 8 is our contribution, mean 
that such examples should/must exist. Note, by the way, that we have been 
unable to provide non-trivial families of MSTD sets which generalise
the sets $A_{2}$ and $A_{10}$.  
\par Otherwise, the outstanding open problem seems to us to be to obtain 
a reasonably tight estimate for 
$L := \limsup_{|A| \rightarrow \infty} f(A)$, where $f$ is the 
function of eq.(1). The set $A_{15}$ leaves the yawning 
gap ${\ln 91 \over \ln 83} \leq L \leq {4 \over 3}$. Finding good 
estimates for the constants $c_{j,k}$ in Theorem 8 is also a very 
appealing problem and several questions in a similar vein are suggested
by the material in [5]. 
\par Finally, one would obviously like to be able to complete the 
classification, up to linear transformation, of MSTD sets of size 9, and
if possible extend the range of computation to larger sizes. An interesting
question is to find the smallest $n_{0} > 0$ 
such that there are infinitely many 
different normalised MSTD sets of size $n_{0}$. All our
constructions, plus those of [3], provide only finitely many sets of any 
given size. An upper bound
is $n_{0} = 16$, for which we can argue as follows : If $A,B \subset 
{\hbox{{\bf Z}}}$ and $C = A \times B \subset {\hbox{{\bf Z}}}^{2}$, then
$|C \pm C| = |A \pm A| \cdot |B \pm B|$. In particular, if
$A$ is an MSTD set and $B$ has at least as many sums as differences, then
$C$ is an MSTD set in {\bf Z}$^{2}$. For suitable $\lambda, \mu 
\in {\hbox{{\bf N}}}$ the map $(a,b) \mapsto \lambda a + \mu b$ takes 
$C$ to a set $C_{\lambda,\mu} \subset {\hbox{{\bf Z}}}$ for which 
$|C_{\lambda,\mu} \pm C_{\lambda,\mu}| = |C\pm C|$ - this is just 
the base expansion method. If $\min\{|A|,|B|\} > 1$ then there is a linear 
transformation taking $C_{\lambda_{1},\mu_{1}}$ to $C_{\lambda_{2},\mu_{2}}$
if and only if $\left| \begin{array}{cc} \lambda_{1} & \mu_{1} \\
\lambda_{2} & \mu_{2} \end{array} \right| = 0$. In particular, we 
get infinitely many different sets up to linear transformation, each of size 
$|A||B|$. Since there is an MSTD set of size 8, namely $A_{1}$, and 
we may take $B = \{0,1\}$ so that $|B+B| = |B-B| = 3$, then we may 
conclude that there are infinitely many normalised 
MSTD sets in {\bf Z} of size 16.         

\begin{center} {\bf References} \end{center}
[1] G.A. Freiman and V.P. Pigaev, {\em The relation between the invariants
$R$ and $T$ (Russian)}, Kalinin. Gos. Univ. Moscow (1973), 172-174.
\newline
[2] J. Marica, {\em On a conjecture of Conway}, Canad. Math. Bull. 
12 (1969), 233-234.
\newline
[3] M.B. Nathanson, {\em Sets with more sums than differences}. Preprint
available online at http://www.arxiv.org/math.NT/0608148
\newline
[4] M.B. Nathanson, {\em Problems in additive number theory, I}. 
Preprint available online at http://www.arxiv.org/math.NT/0604340
\newline
[5] G. Martin and K. O'Bryant, 
{\em Many sets have more sums than differences}. Preprint 
available online at http://www.arxiv.org/math.NT/0608131
\newline
[6] F. Roesler, {\em A mean value density theorem of additive number theory},
Acta Arith. 96 (2000), n0.2, 121-138.
\newline
[7] I.Z. Ruzsa, {\em On the cardinality of $A+A$ and $A-A$}, in : 
Coll. Math. Soc. Bolyai 18 : Combinatorics (Keszthely, 
1976), Akad\'{e}miai Kiad\'{o} (Budapest, 1979), pp. 933-938.
\newline
[8] I.Z. Ruzsa, {\em Sets of sums and differences}, in : 
S\'{e}minaire de Th\'{e}orie des Nombres, Paris 1982/83, Birkh\"{a}user
(1984), pp. 267-273.
\newline
[9] I.Z. Ruzsa, {\em An application of graph theory to additive number theory},
Scientia 3 (1991), 97-109.
\newline
[10] I.Z. Ruzsa, {\em On the number of sums and differences}, 
Acta Math. Hung. 59 (3-4) (1992), 439-447.

\end{document}